\documentstyle{amsppt}
\magnification=1200
\TagsOnRight
\NoBlackBoxes
\topmatter
\NoBlackBoxes
\title On the Boundary Orbit
Accumulation Set\\
for a Domain with Non-Compact\\
Automorphism Group
\endtitle
\rightheadtext{On the Boundary Orbit Accumulation Set}
\footnote[]{{\bf Mathematics
Subject Classification:} 32M05 \hfill}
\footnote[]{{\bf Keywords and Phrases:} Non-Compact 
Automorphism Groups, Boundary Accumulation Points.  \hfill}
\author A. V. Isaev \ \ \ and \ \ \ S. G. Krantz
\endauthor

\abstract For a smoothly bounded pseudoconvex domain $D\subset{\Bbb
C}^n$ of finite type with non-compact holomorphic automorphism group
$\text{Aut}(D)$, we show that the set $S(D)$ of all boundary
accumulation points for $\text{Aut}(D)$ is a compact subset of
$\partial D$ and, if $S(D)$ contains at least three points, it is
a perfect set and thus has the power of the continuum. Moreover, we show that in this case, $S(D)$ is either connected or the number of its connected components is uncountable. We also discuss how $S(D)$ relates to other invariant subsets of $\partial D$.

\endabstract 
\endtopmatter 
\document 
\def\qed{{\hfill{\vrule height7pt width7pt depth0pt}\par\bigskip}}

Suppose that $D\subset{\Bbb C}^n$ is a smoothly bounded domain, 
i.e. $D$ is
bounded and $\partial D$ is $C^\infty$-smooth. We assume that the
group $\text{Aut}(D)$ of holomorphic automorphisms of $D$ is
non-compact; this means, thanks to a classical result of H. Cartan, that
 there is a point $q\in\partial D$ such that for some
$p\in D$ and a sequence $\{f_j\}\subset \text{Aut}(D)$, one has
that $f_j(p)\rightarrow q$ as $j\rightarrow \infty$. Such a point $q$ is
called a {\it boundary accumulation point for}\/ $\text{Aut}(D)$ (see \cite{GK1} for a discussion of this matter).

Let $S(D)$ denote the set of all boundary accumulation points for
$\text{Aut}(D)$. Existing examples of domains with non-compact
automorphism groups (see \cite{FIK} for a discussion of the case of
Reinhardt domains), for which the set $S(D)$ can be found explicitly,
indicate that this set should enjoy some explicit regularity properties. 

For instance, let us for the moment restrict attention to smoothly
bounded Reinhardt domains.  It follows from \cite{FIK} that, for such
a domain $D$, $S(D)$ is always a compact, connected smooth submanifold
of $\partial D$. For such domains one can also observe other
interesting properties of $S(D)$ such as the constancy and minimality
of the rank of the Levi form of $\partial D$ along $S(D)$ (see [H] for
the genesis of these ideas); there is also a certain relation between
this rank and the dimensions of the orbits of the action of
$\text{Aut}(D)$ on $D$. Similarly, the type in the sense of D'Angelo
\cite{D'A1} is constant and maximal along $S(D)$. Many of
these properties, when considered for general domains, appear to be
related to the conjecture of Greene/Krantz \cite{GK2} which states
that every boundary accumulation point for a smoothly bounded domain
must be of finite type.

In the present paper we begin a systematic study of the set $S(D)$ for
a fairly general class of domains, and obtain foundational results on
its topology and the relation of $S(D)$ to other invariant subsets of
$\partial D$. We thank H. Boas, R. Remmert, J. Wolf and S. Fu for stimulating
remarks and suggestions concerning this work. We are also grateful to K. Diederich for a very valuable discussion of the results of this paper.

We say that $\partial D$ is {\it variety-free} at $q\in\partial D$ if
there are no non-trivial germs of complex varieties lying in $\partial
D$ and passing through $q$.

\proclaim{Proposition 1} Let $D$ be a bounded domain in ${\Bbb
C}^n$. Suppose that $\partial D$ is variety-free at each point of
$S(D)$. Then $S(D)$ is compact.  
\endproclaim

\demo{Proof} We need only to prove that $S(D)$ is closed. Let
$\{q_k\}$ be a sequence of points from $S(D)$ such that
$q_k\rightarrow q\in\partial D$ as $k\rightarrow\infty$. 
Since $\partial D$ is variety-free at each point $q_k$, then for every
$q_k$ there is a sequence $\{f_k^j\}$ from $\text{Aut}(D)$ such that
$f_k^j$ converges to the constant map $q_k$ in all of $D$ as
$j\rightarrow\infty$ (see [GK1]).

Fix now a sequence $\{\epsilon_k\}$, $\epsilon_k>0$,
$\epsilon_k\rightarrow 0$ as $k\rightarrow\infty$. Next, fix a point
$p\in D$ and for every $k$ find an index $j(k)$ such that
$|f_k^{j(k)}(p)-q_k|<\epsilon_k$.  It is now obvious that
$f_k^{j(k)}(p)\rightarrow q$ as $k\rightarrow\infty$. Hence $q\in
S(D)$ and $S(D)$ is closed.

The proposition is proved.\qed
\enddemo

\noindent {\bf Remark.} For smoothly bounded domains the 
variety-free assumption
in Proposition 1 would follow from the conjecture of Greene/Krantz.

\proclaim{Theorem 2} Suppose that $D\subset {\Bbb C}^n$ is a smoothly
bounded pseudoconvex domain of finite type. Then, if $S(D)$ contains
at least three points, it is a perfect set and thus has the power of the continuum. Moreover, in this case, $S(D)$ is either connected, or the number of its connected components is uncountable. 
\endproclaim

\demo{Proof}We note that any automorphism of $D$ extends to a
$C^{\infty}$-automorphism of $\overline D$ (see e.g. \cite{D'A2}).

Assume that $S(D)$ contains at least three points.  We will first show 
that $S(D)$ cannot have isolated points. Indeed,
let $q\in S(D)$ be an isolated point. Let $\{f_j\}$ be a sequence in
$\text{Aut}(D)$ such that $f_j\rightarrow q$ in all of $D$ as
$j\rightarrow\infty$. By passing to a subsequence we can also assume
that $f_j^{-1}\rightarrow r$ in all of $D$ where $r\in S(D)$.

Suppose first that $r=q$. Since $S(D)$ contains at least three points,
one can find two distinct points $s,t\in S(D)$, $s,t\ne q$.  Then by
Theorem 1 of \cite{B}, $f_j(s)\rightarrow q$, $f_j(t)\rightarrow q$ as
$j\rightarrow\infty$. Since $q$ is an isolated point of $S(D)$ and
each $f_j$ preserves $S(D)$, we conclude that, for all 
sufficiently large $j$, one has $f_j(s)=f_j(t)=q$. 
This is impossible since every $f_j$ is a
one-to-one mapping on $\overline{D}$.

Suppose now that $r\ne q$. Then there is $s\in S(D)$ such that $s\ne
q,r$. Then by \cite{B}, $f_j(q)\rightarrow q$, $f_j(s)\rightarrow q$
as $j\rightarrow \infty$ which implies as before that for all
sufficiently large $j$, $f_j(q)=f_j(s)=q$, but this is again
impossible since the $f_j$ are one-to-one on $\overline{D}$. 
Thus if $S(D)$ has at least
three elements then $S(D)$ does not have isolated points, and thus is a perfect set.

Assume now that $S(D)$ is disconnected and the number of its connected components is not uncountable. Let $S(D)=\cup_kS_k(D)$ be the decomposition of $S(D)$ into the disjoint union of its connected components. We will show that, for every $k_0$, every $q\
in S_{k_0}(D)$ and every
neighborhood $U$ of $q$ there exists $k_1\ne k_0$ such that $U\cap
S_{k_1}\ne\emptyset$. This implies that the number of connected
components of $S(D)$ is infinite and that each of the sets
$X_m=S(D)\setminus \cup_{k=1}^m S_k(D)$ is open and dense in
$S(D)$. Then since $S(D)$ is compact by Proposition 1 and
since $\cap_{m=1}^\infty X_m=\emptyset$, the Baire Category Theorem gives a
contradiction.

Let $S_{k_0}(D)$ be a component of $S(D)$ and $q\in S_{k_0}$ be such that
there exists a neighborhood $U$ of $q$ that does not contain points
from $S_k(D)$ with $k\ne k_0$. Since $S(D)$ does not have isolated
points, $S_{k_0}(D)$ is not a one-point set. Therefore, by decreasing
$U$ if necessary we can assume that $S_{k_0}(D)\setminus
U\ne\emptyset$. Let $\{f_j\}\subset \text{Aut}(D)$ be such that
$f_j\rightarrow q$, $f_j^{-1}\rightarrow r$ in all of 
$D$, with $r\in S(D)$ as $j\rightarrow\infty$.

Suppose first that $r\in S_{k_0}(D)$. Then by \cite{B}, for any other
connected component $S_{k_1}(D)$ of $S(D)$ with $k_1\ne k_0$, 
one has $f_j(S_{k_1}(D))\subset U$ for all sufficiently large $j$;
this is impossible since $U$ does not contain an entire component of
$S(D)$. If $r\not\in S_{k_0}(D)$, then for all sufficiently large
$j$, $f_j(S_{k_0}(D))\subset U$ which is again impossible.

Thus, $S(D)$ either is connected or has uncountably many components.

The theorem is proved.\qed
\enddemo

As we noted in the proof of Theorem 2 above, for a smoothly bounded
pseudoconvex domain of finite type, the set $S(D)$ is invariant under
the extension of an automorphism of $D$ to the boundary. In the
following proposition we show that $S(D)$ is generically the
smallest invariant subset of $\partial D$.

\proclaim{Proposition 3} Let $D\subset{\Bbb C}^n$ be a smoothly bounded pseudoconvex domain of finite type with non-compact automorphism group. Suppose that $A\subset\partial D$ is non-empty, compact and invariant under $\text{Aut}(D)$. Assume further tha
t $A$ is not a one-point subset of $S(D)$. Then $S(D)\subset A$.

In particular, if $\text{Aut}(D)$ does not have fixed points in $\partial D$, then $S(D)$ is the smallest compact subset of $\partial D$ invariant under $\text{Aut}(D)$.
\endproclaim

\demo{Proof} Since $A$ is closed, it is sufficient to show that every point of $S(D)$ belongs to $\overline{A}$. Let $q\in S(D)$ and $\{f_j\}\subset \text{Aut}(D)$ be such that $f_j\rightarrow q$, $f_j^{-1}\rightarrow r$ in 
all of $D$ as $j\rightarrow \infty$, for some $r\in S(D)$. Since $A$ is not a one-point subset of $S(D)$, there is a point $a\in A$, $a\ne r$. Then, by \cite{B}, $f_j(a)\rightarrow q$ as $j\rightarrow\infty$. Since $A$ is invariant under any $f_j$, we see

 that $f_j(a)\in A$ for all $j$ and thus $q$ is either an accumulation point for $A$, or, if $f_j(a)=q$ for some index $j$, then $q\in A$.

The proposition is proved.\qed
\enddemo

We now derive several corollaries from the above proposition regarding particular\linebreak sets $A$.
 
Fix $0\le k\le n-1$ and denote by $L_k(D)$ the set of all points from $\partial D$ where the rank of the Levi form of $\partial D$ does not exceed $k$. Clearly, each set $L_k(D)$ is a compact subset of $\partial D$ and is invariant under any automorphism 

of $D$. Let $l_1$ denote the minimal rank of the Levi form on $\partial D$ and $l_2$ the minimal rank of the Levi form on $\partial D\setminus L_{l_1}(D)$. For these sets, Proposition 3 gives the following corollary (that was first proved in \cite{H}). No
te that the proof in \cite{H} was also based on the results of \cite{B}.
\pagebreak

\proclaim{Corollary 4} Let $D$ be as in Proposition 3. Then either
\smallskip 

\noindent (i) $S(D)\subset L_{l_1}(D)$,

\noindent or

\noindent (ii) $L_{l_1}(D)$ is a one-point subset of 
$S(D)$ and $S(D)\subset L_{l_2}(D)$.
\endproclaim

\demo{Proof} If $L_{l_1}(D)$ is not a one-point subset of $S(D)$ then,
by Proposition 3, $S(D)\subset L_{l_1}(D)$. Suppose now that
$L_{l_1}(D)$ is a one-point subset of $S(D)$. Then, since $L_{l_1}(D)$
is strictly contained in $L_{l_2}(D)$, one has $S(D)\subset
L_{l_2}(D)$.

The corollary is proved.\qed
\enddemo

By a similar argument, one can endeavor to prove an analogous property of
the type $\tau(q)$, $q\in\partial D$, in the sense of
D'Angelo. Indeed, denote by $T_k(D)$ the set of all points
$q\in\partial D$ where $\tau(q)$ is at least $k$. We choose $t_1$ and
$t_2$ such that $T_{t_1}(D)\ne\emptyset$, $t_2<t_1$, and there exists
a point of type $t_2$ in $\partial D\setminus T_{t_1}(D)$. Since
$\tau$ is invariant under automorphisms of $D$, so is every set
$T_k(D)$. However, the sets $T_k(D)$ do not have to be closed, as the
type function $\tau$ may not be upper-semicontinuous on $\partial D$
(see e.g. an example in \cite{D'A2}, p. 136). Therefore, for the type
we only have a somewhat weaker result.

\proclaim{Corollary 5} Let $D$ be as in Proposition 3. Then either
\smallskip 

\noindent (i) $S(D)\subset \overline{T_{t_1}(D)}$,

\noindent or

\noindent (ii) $T_{t_1}(D)$ is a one-point subset of $S(D)$ and $S(D)\subset \overline{T_{t_2}(D)}$.
\endproclaim

In place of the type function $\tau$, one can consider the
multiplicity function $\mu$ on $\partial D$ (see \cite{D'A2}, p. 145
for the definition), which is also invariant under the extensions of
automorphisms to $\partial D$. It should be noted that, for
$q\in\partial D$, the number $\tau(q)$ is finite if and only if
$\mu(q)$ is finite. In contrast with $\tau$, however, the function
$\mu$ is upper-semicintinuous on $\partial D$. Analogously to what we
have done above for the function $\tau$, denote by $M_k(D)$ the set of
all points $q\in\partial D$, where $\mu(q)$ is at least $k$ and choose
$m_1$ and $m_2$ such that $m_1=\text{max}_{q\in\partial D}\,\mu(q)$,
$m_2<m_1$, and there exists a point of multiplicity $m_2$ in $\partial
D\setminus M_{m_1}(D)$. Due to the upper-semicintinuity and invariance
of $\mu$, each set $M_k(D)$ is a compact subset of $\partial D$
that is invariant under $\text{Aut}(D)$. This observation gives 
the following analogue of Corollary 6 for $M_{m_1}, M_{m_2}$.

\proclaim{Corollary 6} Let $D$ be as in Proposition 3. Then either
\smallskip 

\noindent (i) $S(D)\subset M_{m_1}(D)$,

\noindent or

\noindent (ii) $M_{m_1}(D)$ is a one-point subset of $S(D)$ and $S(D)\subset M_{m_2}(D)$.
\endproclaim

The proof of Corollaries 5, 6 is completely analogous to that of Corollary 4.
\medskip

\noindent {\bf Remarks.}
\medskip

1. It is plausible that Theorem 2 and Corollaries 4--6 hold without
the assumptions of pseudoconvexity and finite type.

2. We note that, in complex dimension 2, the type $\tau$ is
upper-semicontinuous. As a result, Corollary 5 can be stated in this
case without passing to the closures of the $T_{t_j}$. Also, in complex
dimension 2, Corollary 5 is a consequence of the explicit
classification of smoothly bounded pseudoconvex domains of finite type
with non-compact automorphism groups \cite{BP1}.

3. For a smoothly bounded circular domain, the set $S(D)$ clearly
cannot be a one- or two-point set. Thus, Theorem 2 gives that, for smoothly bounded pseudoconvex circular domains of finite type, $S(D)$ is always a perfect set. Next, since for such domains the automorphism group cannnot have fixed points on the boundary
, Corollary 4 implies that in this case the Levi form of $\partial D$ has constant rank along $S(D)$ and minimizes its rank over $\partial D$ on $S(D)$ (see also \cite{H}). It also should be noted here that, by the results of \cite{BP2}, every smoothly bo
unded convex domain of finite type with non-compact automorphism group is biholomorphically equivalent to a certain polynomially defined domain that admits an action of the
two-dimensional torus ${\Bbb T}^2$. Therefore, for any such a domain,
$S(D)$ also is a perfect set, and the rank of the Levi form is constant and minimal on $S(D)$.
 
4. The results of \cite{FIK} imply that, for a smoothly bounded
Reinhardt domain $D$, the type is constant along $S(D)$ and maximizes
on $S(D)$ the type over $\partial D$. It is an interesting question
whether there exists an analogue of this fact for more general domains
(cf. Corollaries 5, 6). Note that one can make a statement analogous to Corollary 6 for the multitype introduced in \cite{C}, since the multitype function is upper-semicontinuous with respect to lexicographic ordering.

5. It also follows from \cite{FIK} that, for a smoothly bounded
Reinhardt domain $D$, the real dimension of any orbit of the action of
$\text{Aut}(D)$ on $D$ is at least $2(k+1)$, where $k$ is the rank of
the Levi form of $\partial D$ along $S(D)$. Moreover, there is
precisely one orbit of minimal dimension $2(k+1)$ (see \cite{K} for a
discussion of this phenomenon). Also, the orbit of minimal dimension
approaches every point of $S(D)$ non-tangentially, whereas any other
orbit approaches every point of $S(D)$ only along tangential
directions. It would be interesting to know if similar statements hold
for more general, e.g. circular, domains. The fact that there exists
an orbit that approaches $S(D)$ non-tangentially would be very
important for a proof of the Greene/Krantz conjecture. It also could
be used to show that $S(D)$ is a smooth submanifold of $\partial D$.
\medskip

We wish to conclude this paper with a list of immediate open problems
that arise from the above discussion and which complement some of the
preceding remarks.  
\medskip

{\bf Open Problems.}
\medskip

1. For a smoothly ($C^\infty$) bounded domain $D\subset {\Bbb C}^n$, 
can the set $S(D)$ be a one- or two-point set? Note that the 
reference [GK3] gives an example of a domain with $C^{1-\epsilon}$ boundary, for which $S(D)$ has only two points.  It appears that this example
can be modified, using a parabolic group of automorphisms, so that
$S(D)$ has just one point.  We shall explore this matter further, and
additionally investigate increasing the boundary smoothness, in
a future paper. Indications are that the case of finite boundary
smoothness will be different from the case of infinite boundary
smoothness.

2. For a smoothly bounded domain $D$, can the set $S(D)$ have uncountably many components, for example, can it be a Cantor-type set?

3. For a smoothly bounded domain $D$, is the set $S(D)$ always a
smooth submanifold of $\partial D$? Note that the results of
\cite{FIK} imply that, for a smoothly bounded Reinhardt domain, $S(D)$
is always a smooth submanifold of $\partial D$ that is diffeomorphic to a
sphere of odd dimension.

4. Is it always true that the rank of the Levi form is in fact
constant and minimal along $S(D)$ and that the type is constant and
maximal along $S(D)$?  
\medskip

This work was completed while the first 
author was an Alexander von Humboldt Fellow at the 
University of Wuppertal.  Research at MSRI by the second author was
supported in part by NSF Grant DMS-9022140.

\bigskip

{\obeylines
Centre for Mathematics and Its Applications 
The Australian National University 
Canberra, ACT 0200
AUSTRALIA 
E-mail address: Alexander.Isaev\@anu.edu.au
\smallskip
and

Bergische Universit\"at
Gesamthochschule Wuppertal
Mathematik (FB 07)
Gaussstrasse 20
42097 Wuppertal
GERMANY
E-mail address: Alexander.Isaev\@math.uni-wuppertal.de
\bigskip

Department of Mathematics
Washington University, St.Louis, MO 63130
USA 
E-mail address: sk\@math.wustl.edu
\smallskip
and
\smallskip
MSRI
1000 Centennial Drive
Berkeley, CA 94720
USA
E-mail address: krantz\@msri.org}

\enddocument